\documentclass[12pt]{amsart}

\usepackage{times,amsfonts,amsmath,amstext,amsbsy,amssymb,
amsopn,amsthm,upref,eucal, mathrsfs}

\usepackage{url}

\def\thispartnum{II}
\input extref-arxiv.tex

\newcommand \buf {{\mathfrak {B}}}

\newcommand \bbuf {{\tilde {\mathfrak {B}}}}

\DeclareMathOperator{\Int}{int}

\DeclareMathOperator{\tr}{tr}

\newcommand \Prob {{\mathbb P}}

\DeclareMathOperator{\Leb}{Leb}

\newcommand  \loc {{\mathrm{loc}}}

\newcommand  \fin {{\mathrm{fin}}}
\newcommand  \Mfin {\mathfrak{M}_\fin}

\DeclareMathOperator{\ang}{angle}

\newcommand \Conf {{\mathrm {Conf}}}

\newcommand\scrI{{\mathscr I}}

\newtheorem{theorem}{Theorem}[section]

\newtheorem{corollary}[theorem]{Corollary}
\newtheorem{proposition}[theorem]{Proposition}

\begin{document}
\title[The Ergodic Decomposition of Infinite Pickrell Measures. II]{Infinite Determinantal Measures and
The Ergodic Decomposition of Infinite Pickrell Measures II. Convergence of determinantal measures}

\author{Alexander I. Bufetov}

\date{}

\address{ Aix-Marseille Universit{\'e}, Centrale Marseille, CNRS,   I2M}
\address{UMR 7373, rue F. Joliot Curie , Marseille, France}
\address{Steklov Institute of Mathematics,
Moscow, Russia}
\address{Institute for Information Transmission Problems,
 Moscow, Russia}
\address{National Research University Higher School of Economics,
 Moscow, }
 \address{Russia}

\begin{abstract}
The second part of the paper mainly deals with convergence of infinite 
determinantal measures, understood as the convergence of the approximating finite determinantal measures.
In addition to the usual weak topology on the space of probability measures on the space of configurations, we also consider 
the weak topology on the space of finite measures on the space of finite measures on the half-line, used via 
the natural  immersion, well-defined almost surely with respect to the infinite Bessel point process, 
of the space of configurations into the space of finite measures on the half-line. The main results of the second part are sufficient conditions for tightness of families of determinantal measures, for convergence of sequences of induced proceses, 
as well as for sequences of finite-dimensional perturbations of determinantal processes. 
\end{abstract}

\maketitle
\tableofcontents
\section{Introduction}

This paper is the second of the cycle of three papers giving the explicit construction of the ergodic decomposition of infinite Pickrell measures.
Quotes to the other parts of the paper \cite{infdet1, infdet3} are organized as follows: Corollary I.1.10, equation (III.15), etc.

The paper is organized as follows.
Section 2 studies convergence of determinantal probability measures
given by positive contractions that are locally trace-class.
We start by recalling that locally trace-class convergence of operators
implies weak convergence of the corresponding determinantal measures
in the space of probability measures on the space of configurations.
In the study of infinite Pickrell measures, we need to consider
induced processes  of the Bessel point process
as well as  as finite-rank perturbations of the Bessel point process, and
in Section 2 sufficient conditions  are given for the convergence of induced processes and of processes induced by finite-rank perturbations. We conclude Section 2 by
establishing, for infinite determinantal measures obtained as finite-rank perturbations,
the convergence of the family of  determinantal processes obtained by inducing on an exhausting family of
subsets of the phase space to the initial, unperturbed, determinantal process.

In Section 3, we embed suitable subsets of the space of configurations
into the space of {\it finite} measures on the phase space $E$ and
give sufficient conditions for precompactness of families of determinantal
measures with respect to the weak topology on the space of finite measures on the space of {\it finite} measures on $E$
(which is stronger than the usual weak topology on the space of finite measures
on the space of {\it Radon} measures, equivalent to the weak topology on the space of finite measures
on the space of configurations).
This  step is needed  for proving the vanishing of
the ``Gaussian parameter'' for the ergodic components of Pickrell measures.
 Borodin and Olshanski \cite{BO} proved  this vanishing for the ergodic components of Hua-Pickrell measures: in fact, the estimate of their argument can be interpreted as the assertion of {\it tightness} of the family of rescaled radial parts of Hua-Pickrell measures considered as measures in the space of finite measures on the space
of {\it finite} measures. We next study weak convergence of induced processes and of finite-rank perturbations with respect to the new topology.

\section{Convergence of determinantal measures}

\subsection{Convergence of operators and convergence of measures}

We consider determinantal probability measures induced by positive contractions and start by  recalling that  convergence of a sequence of such operators in the
space of locally trace-class operators implies the
weak convergence of corresponding determinantal probability measures in the space of finite measures on the space of configurations.

\begin{proposition} \label{optomeas} Assume that the operators $K_n\in \scrI_{1,\loc}\bigl(E,\mu\bigr)$, $n\in{\mathbb N}$, $K\in \scrI_{1,\loc}\bigl(E,\mu\bigr)$ induce determinantal probability measures
$\Prob_{K_n}$, $n\in{\mathbb N}$, $\Prob_K$  on $\Conf(E)$.
If $K_n\to K$  in $\scrI_{1,\loc}\bigl(E,\mu\bigr)$ as $n\to\infty$,
then $\Prob_{K_n}\to\Prob_K$ with respect to the weak topology on $\Mfin(\Conf(E))$ as $n\to\infty$.
\end{proposition}
This proposition is immediate from the definition of determinantal probability measures and Proposition \ref{fin-distr-unique}.
From the classical Heine-Mehler asymptotics (cf. Proposition \ref{tilde-kernel-kns-js}) we now have the following immediate
\begin{corollary}
For any $s>-1$, we  have
$$
{\tilde K}_n^{(s)}\to {\tilde J}_{s}\text{ in }\scrI_{1, \loc}((0, +\infty), \Leb)
$$
and
$$
\Prob_{{\tilde K}_n^{(s)}}\to \Prob_{{\tilde J}_{s}}\text{ in }\Mfin\Conf((0, +\infty)).
$$
\end{corollary}

Our next aim is to show that, under certain additional assumptions, the convergence above persists under passage to induced processes as well as to finite-rank perturbations. We proceed to precise statements.
\subsection{Convergence of induced processes}
Recall that if $\Pi$ is a projection operator acting on $L_2(E, \mu)$ and $g$ is a nonnegative bounded measurable function on $E$ such that the operator
$1+(g-1)\Pi$ is invertible, then we have set
$$
\tilde{\buf}(g, \Pi)=\sqrt{g}\Pi(1+(g-1)\Pi)^{-1}\sqrt{g}.
$$
We now fix $g$ and  establish the connection between convergence of the sequence  $\Pi_n$ and the corresponding
sequence  $\tilde{\buf}(g, \Pi_n)$.

\begin{proposition} \label{conv-induced} Let $\Pi_n,\Pi\in\scrI _ {1,\loc}(E,\mu)$ be orthogonal projection operators,
and let $g: E\to [0,1]$ be a measurable function such that
$$
 \sqrt{1-g}\Pi\sqrt{1-g}\in\scrI _ {1}(E, \mu), \sqrt{1-g}\Pi_n\sqrt{1-g}\in\scrI _ {1}(E, \mu), n\in {\mathbb N}.
$$
Assume furthermore that
\begin{enumerate}
\item $\Pi_n\to\Pi$ in
$\scrI_{1,\loc}(E, \mu)$ as $n\to\infty$;
\item $\lim\limits_{n\to\infty}\tr\sqrt{1-g}\Pi_n\sqrt{1-g}=
\tr\sqrt{1-g}\Pi\sqrt{1-g}$;
\item  the operator $1+(g-1)\Pi$ is invertible.
\end{enumerate}
Then the operators $1+(g-1)\Pi_n$ are also invertible for all sufficiently large $n$, and we have
$$\tilde{\buf}(g, \Pi_n)\to \tilde{\buf}(g, \Pi)\text{ in }\scrI_{1,\loc}(E,\mu)$$
and, consequently,
$$\Prob_{\tilde{\buf}(g, \Pi_n)}\to\Prob_{\tilde{\buf}(g, \Pi)}$$
with respect to the weak topology on $\Mfin(\Conf(E))$ as $n\to\infty$.

\end{proposition}

{\bf Remark.} The second requirement could have been replaced by the requirement that $(g-1)\Pi_n$ converge to $(g-1)\Pi$ in norm,
which is weaker  and is what we shall actually use; nonetheless, in applications it will be more convenient to check the convergence of traces rather than the norm convergence of operators.

Proof.  The first two requirements and Gr{\"u}mm's
Theorem (see Simon \cite{Simon}) imply that
$$\sqrt{1-g}\Pi_n\to \sqrt{1-g}\Pi\text{ in }\scrI_2(E, \mu),$$ whence, a fortiori,   $$(g-1)\Pi_n\to (g-1)\Pi$$
in norm as $n\to\infty$.
We now take a bounded Borel subset $D\subset E$
and check that, as $n\to\infty$, we have
\begin{equation}\label{bg-conv}
\chi_D\bbuf(g, \Pi_n)\chi_D\rightarrow\chi_D\tilde{\buf}(g, \Pi)\chi_D\text{ in }\scrI_{1}(E,\mu).
\end{equation}

Our assumptions directly imply the norm convergence
\begin{equation}\label{inv-conv}
(1+(g-1)\Pi_n)^{-1}\to(1+(g-1)\Pi)^{-1}.
\end{equation}
Furthermore, $$\chi_D\Pi_n\to\chi_D\Pi$$ as  $n\to\infty$ in the
strong operator topology; besides, by our assumptions, we have
$$\lim\limits_{n\to\infty}
\tr\chi_D\Pi_n\chi_D=\tr\chi_D\Pi\chi_D,$$ whence, by Gr{\"u}mm's
Theorem ,  we have $\chi_D\Pi_n\to\chi_D\Pi$ in Hilbert-Schmidt norm,
 and, a fortiori,  in norm.

It follows that the convergence (\ref{bg-conv}) also takes place in norm. To verify the desired $\scrI_1$ convergence, by
Gr{\"u}mm's Theorem again, it suffices to check the relation
\begin{equation}
\label{tr-bg-conv}
\lim\limits_{n\to\infty}\tr\chi_D\bbuf(g,\Pi_n)\chi_D=\tr\chi_D\bbuf(g, \Pi)\chi_D.
\end{equation}

First,
if $A$ is a bounded operator, and  $K_1, K_2\in \scrI_2$,  then one directly verifies the inequality
$$
\left|\tr(K_1^*AK_2)\right| \leq||K_1||_{\scrI_2}\cdot  ||A|| \cdot||K_2||_{\scrI_2}.
$$

It easily  follows that the function $\tr(K_1^*AK_2)$ is continuous
as long as $K_1, K_2$ are operators in $\scrI_2$, and $A$ is  a bounded operator.
The desired convergence of traces (\ref{tr-bg-conv})
follows from the said continuity since
$$\chi_D\bbuf(g,\Pi)\chi_D=\chi_D\sqrt{g-1}\Pi(1+(g-1)\Pi)^{-1}\Pi\sqrt{g-1}\chi_D,
$$
and we have the norm convergence (\ref{inv-conv}) and the
$\scrI_2$-convergence
$$
\chi_D\Pi_n\to \chi_D\Pi.
$$

%We first check that the convergence (\ref{bg-conv})
%takes place in norm.

\subsubsection{Convergence of finite-rank perturbations}
We now proceed to the study of convergence of finite-rank perturbations of locally trace-class projection operators.
Let $L_n$, $L\subset L_2(E,\mu)$ be closed subspaces, and let $\Pi_n,\, \Pi$ be the corresponding orthogonal projection operators. Assume we are given non-zero vectors $v^{(n)}\in L_2(E,\mu)$, $n\in\mathbb{N}$, $v\in L_2(E,\mu)$, and let $\widetilde{\Pi}_n, \widetilde{\Pi}$ be the operators of orthogonal projection onto, respectively, the subspaces $L_n+\mathbb{C}v^{(n)}$, $n\in\mathbb{N}$ and $L\oplus\mathbb{C}v$ .

\begin{proposition}\label{one-pert}
Assume
\begin{enumerate}
\item $\Pi_n\to \Pi$ in the strong operator topology as $n\to\infty$;
\item $v^{(n)}\to v$ in $L_2(E,\mu)$ as $n\to\infty$;
\item $v\notin L$.
\end{enumerate}
Then ${\widetilde \Pi}_n\to {\widetilde \Pi}$ in the strong operator topology as $n\to\infty$.

If, additionally,
$$
\Pi_n\to \Pi \  \text{in} \ \scrI_{1,\loc}(E,\mu) \ \text{as}  \ n\to\infty,
$$
then also
$$
\widetilde{\Pi}_n\to \widetilde{\Pi} \ \text{in} \ \scrI_{1,\loc}(E,\mu) \ \text{as}  \ n\to\infty.
$$

\end{proposition}
%\end{proposition}

Let $\ang(v, H)$ stands for the angle between a vector $v$ and a subspace $H$.
Our assumptions imply that there exists $\alpha_0>0$ such that for sufficiently large $n$
%добавлено в англ.версию: "for sufficiently large $n$"
$$
\ang(v_n, L_n)\geq \alpha_0.
$$

Decompose
$$v^{(n)}=\beta(n)\widetilde{v}^{(n)}+\widehat{v}^{(n)},$$
where $\widetilde{v}^{(n)}\in L_n^{\perp}$,\, $\|\widetilde{v}^{(n)}\|=1$, $\widehat{v}^{(n)}\in L_n$. In this case we have
$$\widetilde{\Pi}_n=\Pi_n+P_{\widetilde{v}^{(n)}}\,,$$
where $P_{\widetilde{v}^{(n)}}\colon v\to\langle v,\widetilde{v}^{(n)}\rangle\widetilde{v}^{(n)}\,,$ is the operator of the orthogonal projection onto the subspace $\mathbb{C}\widetilde{v}^{(n)}$.

Similarly, decompose
$$v=\beta\widetilde{v}+\widehat{v}$$
with $\widetilde{v}\in L^{\perp}$,\, $\|\widetilde{v}\|=1$\,, $\widehat{v}\in L$, and, again, write
$$\widetilde{\Pi}_n=\Pi_n+P_{\widetilde{v}},$$
with $P_{\widetilde{v}}(v)=\langle v,\widetilde{v}\rangle\widetilde{v}$.

Our assumptions 2 and 3 imply that $\widetilde{v}^{(n)}\to\widetilde{v}$ in $L_2(E,\mu)$. It follows that $P_{\widetilde{v}^{(n)}}\to P_{\widetilde{v}}$ in the strong operator topology and
also, since our operators have one-dimensional range, in $\scrI_{1,\loc}(E,\mu)$, which implies the proposition.

The case of perturbations of higher rank follows by induction.
Let $m\in\mathbb{N}$ be arbitrary and assume we are given non-zero vectors $v_1^{(n)}, v_2^{(n)}, \dots, v_m^{(n)}\in L_2(E,\mu)$, $n\in\mathbb{N}$, $v_1, v_2, \dots, v_m\in L_2(E,\mu)$.
Let $$\widetilde{L}_n=L_n\oplus\mathbb{C}v_1^{(n)}\oplus\dots\oplus\mathbb{C}v_m^{(n)},$$
$$\widetilde{L}=L\oplus\mathbb{C}v_1\oplus\dots\oplus\mathbb{C}v_m,$$
and let $\widetilde{\Pi}_n, \widetilde{\Pi}$ be the corresponding projection operators.

Applying Proposition \ref{one-pert} inductively, we obtain
\begin{proposition} \label{fin-pert} Assume
\begin{enumerate}
\item $\Pi_n\to \Pi$ in the strong operator topology as $n\to\infty$;
\item $v_i^{(n)}\to v_i$ in $L_2(E,\mu)$ as $n\to\infty$ for any $i=1,\dots,m$\,;
\item $v_k\notin L\oplus\mathbb{C}v_1\oplus\dots\oplus\mathbb{C}v_{k-1}$, $k=1, \dots, m$.
\end{enumerate}
Then ${\widetilde \Pi}_n\to {\widetilde \Pi}$ in the strong operator topology as $n\to\infty$.
If, additionally,
$$
\Pi_n\to \Pi \  \text{in} \ \scrI_{1,\loc}(E,\mu) \ \text{as}  \ n\to\infty,
$$
then also
$$
\widetilde{\Pi}_n\to \widetilde{\Pi} \ \text{in} \ \scrI_{1,\loc}(E,\mu) \ \text{as}  \ n\to\infty,
$$
and, consequently,
$\Prob_{\widetilde{\Pi}_n}\to\Prob_{\widetilde{\Pi}}$ with respect to the weak topology on $\Mfin(\Conf(E))$ as $n\to\infty$.
\end{proposition}

\subsection{Application to infinite determinantal measures}

Take a sequence $$\mathbb{B}^{(n)}=\mathbb{B}\left(H^{(n)},E_0\right)$$ of
infinite determinantal measures with $H^{(n)}=L^{(n)}+V^{(n)}$, where
$L^{(n)}$ is, as before, the range of a projection operator
$\Pi^{(n)}\in\scrI_{1,\loc}(E,\mu)$, and $V^{(n)}$ is
finite-dimensional. Note that the subset $E_0$ is fixed throughout.

Our aim is to give sufficient conditions for convergence of
$\mathbb{B}^{(n)}$ to a limit measure
$\mathbb{B}=\mathbb{B}\left(H,E_0\right)$, $H=L+V$, the subspace $L$ being
the range of a projection operator $\Pi\in\scrI_{1,\loc}(E,\mu)$.

\begin{proposition}\label{conv-inf-det}  Assume
\begin{enumerate}
\item $\Pi^{(n)}\to\Pi$ in $\scrI_{1,\loc}(E,\mu)$ as
$n\to\infty$\,;
\item the subspace $V^{(n)}$ admits a basis $v_1^{(n)},\dots,v_m^{(n)}$
and the subspace $V$ admits a basis $v_1,\dots,v_m$ such that
$$v_i^{(n)}\to v_i\;\text{ in }\; L_{2, \loc}(E,\mu)\;\text{ as }\;
n\to\infty \;\text{ for all }\; i=1,\dots,m\;.$$
\end{enumerate}

Let $g\colon E\to [0,1]$ be a positive measurable function such that
\begin{enumerate}
\item
$\sqrt{1-g}\Pi^{(n)}\sqrt{1-g}\in\scrI_1(E,\mu)$,\;
$\sqrt{1-g}\Pi\sqrt{1-g}\in\scrI_1(E,\mu)$\,;
\item
$\lim\limits_{n\to\infty}\tr\sqrt{1-g}\Pi^{(n)}\sqrt{1-g}=\tr\sqrt{1-g}\Pi\sqrt{1-g}$\,;
\item $\sqrt{g}V^{(n)}\subset L_2(E,\mu)$,\; $\sqrt{g}V\subset
L_2(E,\mu)$\,;
\item $\sqrt{g}v_i^{(n)}\to\sqrt{g}v_i$ in $L_2(E,\mu)$ as
$n\to\infty$ for all $i=1,\dots,m$\,.
\end{enumerate}

Then

\begin{enumerate}
\item the subspaces $\sqrt{g}H^{(n)}$ and $\sqrt{g}H$ are closed\,;
\item the operators $\Pi^{(g,n)}$ of orthogonal projection onto the
subspace $\sqrt{g}H^{(n)}$ and the operator $\Pi^g$ of orthogonal
projection onto the subspace $\sqrt{g}H$ satisfy
$$\Pi^{(g,n)}\to\Pi^g \;\text{ in }\; \scrI_{1,\loc}(E,\mu)
\;\text{ as }\; n\to\infty\;.$$
\end{enumerate}

\end{proposition}

\begin{corollary} In the notation and under the assumptions of  Proposition \ref{conv-inf-det}, we have
\begin{enumerate}
\item $\Psi_g\in L_1(\Conf(E), \mathbb{B}^{(n)})$ for all $n$, $\Psi_g\in L_1(\Conf(E), \mathbb{B})$;

\item $$
\frac{\Psi_g\mathbb{B}^{(n)}}{\displaystyle \int\limits_{\Conf(E)} \Psi_g\,d\mathbb{B}^{(n)}}\to \frac{\Psi_g\mathbb{B}}{\displaystyle \int\limits_{\Conf(E)}\Psi_g\,d\mathbb{B}}
$$ with respect to the weak topology on $\Mfin(\Conf(E))$ as $n\to\infty$.
\end{enumerate}
\end{corollary}

Indeed, the Proposition and the Corollary are immediate from the characterization of multiplicative
functionals of infinite determinantal  measures given in Proposition \ref{mult-he1} and Corollary \ref{fin-rank-mult}, the  sufficient conditions of convergence of induced processes and finite-rank perturbations
given in Propositions \ref{conv-induced}, \ref{fin-pert}, and the characterization of convergence with respect to the weak topology on $\Mfin(\Conf(E))$ given in Proposition \ref{optomeas}.

\subsection{Convergence of approximating kernels and the proof of Proposition \ref{asympt-R-bessel}}

Our next aim is to show that, under certain additional assumptions, if a sequence $g_n$ of measurable functions converges to $1$, then the operators $\Pi^{g_n}$ defined by \eqref{def-Pi-g} converge to $Q$ in $\scrI_{1, \loc}(E, \mu)$.

Given  two closed subspaces $H_1, H_2$ in $L_2(E, \mu)$, let $\ang(H_1, H_2)$
be the angle between $H_1$ and $H_2$, defined as the infimum of angles between all nonzero vectors
in $H_1$ and $H_2$; recall that if one of the subspaces has finite dimension, then the infimum is achieved.
\begin{proposition} \label{gn-conv}
Let $L$, $V$, and $E_0$ satisfy Assumption~\ref{lve}, and assume additionally that we have $V\cap L_2(E,\mu)=0.$
Let $g_n:E\to (0,1]$ be a sequence of positive measurable functions such that
\begin{enumerate}
\item for all $n\in {\mathbb N}$ we have
$\sqrt{1-g_n}Q\sqrt{1-g_n}\in\scrI_1(E,\mu)$;\item for all $n\in {\mathbb N}$ we have $\sqrt{g_n}V\subset L_2(E, \mu)$;
\item there exists $\alpha_0>0$ such that for all $n$ we have
$$
\ang(\sqrt{g_n}H, \sqrt{g_n}V)\geq \alpha_0;
$$
\item for any bounded $B\subset E$  we have
$$
\inf\limits_{n\in {\mathbb N}, x\in E_0\cup B} g_n(x)>0;\
\lim\limits_{n\to\infty} \sup\limits_{x\in E_0\cup B} \left|g_n(x)-1\right|=0.
$$
\end{enumerate}
Then, as $n\to\infty$, we have
$$
\Pi^{g_n}\to Q\text{ in }\scrI_{1, \loc}(E, \mu).
$$
\end{proposition}

Using the second remark after Theorem \ref{infdet-he}, one can extend Proposition \ref{gn-conv} also to nonnegative functions that admit zero values. Here we restrict ourselves to characteristic functions of the form $\chi_{E_0\cup B}$ with $B$ bounded, in which case we have the following

\begin{corollary}
\label{cor-conv}
Let $B_n$ be an increasing sequence of bounded Borel sets exhausting $E\setminus E_0$.
 If there exists $\alpha_0>0$ such that for all $n$ we have
$$
\ang(\chi_{E_0\cup B_n}H, \chi_{E_0\cup B_n}V)\geq \alpha_0,
$$
then
$$
\Pi^{E_0\cup B_n}\to Q\text{ in }\scrI_{1,\loc}(E, \mu).
$$
\end{corollary}

Informally, Corollary \ref{cor-conv} means that, as $n$ grows, the
induced processes of our determinantal measure
on subsets $\Conf(E; E_0\cup B_n)$ converge to the ``unperturbed''
determinantal point process $\Prob_{Q}$.

Note that Proposition \ref{asympt-R-bessel} is an immediate corollary
of Proposition \ref{gn-conv} and Corollary \ref{cor-conv}.

Proof of Proposition \ref{gn-conv}.

We start by showing that, as $n\to\infty$, we have $g_nQ\to Q$ in norm.

Indeed, take $\varepsilon>0$ and choose a bounded set $B_{\varepsilon}$ in
such a way that
$$
\tr\chi_{E\backslash(E_0\cup B_{\varepsilon})}Q\chi_{E\backslash(E_0\cup B_{\varepsilon})}<\frac{\varepsilon^2}4.
$$
Since $g_n\to1$ uniformly on $E_0\cup B_{\varepsilon}$, we have
$$\chi_{E_0\cup B_{\varepsilon}}(g_n-1)Q\to0$$ in norm as $n\to\infty$.
Furthermore, we have
$$
\|\chi_{E\backslash(E_0\cup B_{\varepsilon})}Q\| \leq \|\chi_{E\backslash(E_0\cup B_{\varepsilon})}Q\|_{\scrI_2}<\frac{\varepsilon}2\,.
$$

Consequently, for $n$ sufficiently big, we have:
$$\|(g_n-1)Q\|\leq\|\chi_{E_0\cup B_{\varepsilon}}(g_n-1)Q\|+\|\chi_{E\backslash(E_0\cup
B_{\varepsilon})}Q\|<\varepsilon\,,$$
and, since $\varepsilon$ is arbitrary, we have, as desired, that $g_nQ\to
Q$ in norm as $n\to\infty$.

In particular, we have
$$\left(1+(g_n-1)Q\right)^{-1}\to1$$
in norm as $n\to\infty$.

Now, since $g_n\to1$ uniformly on bounded sets, for any bounded Borel
subset $B\subset E$, we have
$$\chi_B\sqrt{g_n}Q\to\chi_BQ\;\;\,\text{in}\;\;\scrI_2(E,\mu)$$
as $n\to\infty$. Consequently, we have
$$\chi_B\sqrt{g_n}Q\left(1+(g_n-1)Q\right)^{-1}Q\sqrt{g_n}\chi_B\to\chi_BQ\chi_B$$
in $\scrI_1(E,\mu)$ as $n\to\infty$, and, since $B$ is arbitrary, we
obtain
$$Q^{g_n}\to Q\;\;\,\text{in}\;\;\scrI_{1,\loc}(E,\mu)\,.$$

We now let $V_n$ be the orthogonal complement of $\sqrt{g_n}L$ in
$\sqrt{g_n}L+\sqrt{g_n}V$, and let $\tilde{P}^{(n)}$ be the operator of
orthogonal
projection onto $V_n$.

By definition, we have
$$\Pi^{g_n}=Q^{g_n}+\tilde{P}^{(n)}\,.$$

To complete the proof, it suffices to establish that, as $n\to\infty$, we
have
$$\tilde{P}^{(n)}\to0\;\;\,\text{in}\;\;\scrI_{1,\loc}(E,\mu)\,,$$
to do which, since $\tilde{P}^{(n)}$ are projections onto subspaces whose
dimension does not exceed that of $V$, it suffices to show that for any
bounded set $B$
we have $\tilde{P}^{(n)}\to 0$ in strong operator topology as $n\to\infty$.

Since the angles between subspaces $\sqrt{g_n}L$ and $\sqrt{g_n}V$ are
uniformly bounded from below, it suffices to establish the strong
convergence to 0
of the operators $P^{(n)}$ of orthogonal projections onto the
subspaces $\sqrt{g_n}V$.

Let, therefore, $\varphi\in L_2(E,\mu)$ be supported in a bounded Borel set
$B$; it suffices to show that $\|P^{(n)}\varphi\|\to0$ as $n\to\infty$.
But since $V\cap L_2(E,\mu)=0$, for any $\varepsilon>0$ there exists a
bounded set $B_{\varepsilon}\supset B$ such that for any $\psi\in V$
we have
$$\displaystyle\frac{\|\chi_B\psi\|}{\|\chi_{B_{\varepsilon}}\psi\|}<\varepsilon^2.$$

 We have
\begin{multline*}
\|\Pi^{B_{\varepsilon}}\varphi\|^2=\langle\varphi,\Pi^{B_{\varepsilon}}\varphi\rangle=\\
=\langle\varphi,\chi_B\Pi^{B_{\varepsilon}}\varphi\rangle\leq\|\varphi\|\,\|\chi_B\Pi^{B_{\varepsilon}}\varphi\|\leq\\
\leq\|\varphi\|\,\varepsilon\|\Pi^{B_{\varepsilon}}\varphi\|\leq\varepsilon\|\varphi\|\,\|\Pi^{B_{\varepsilon}}\varphi\|\leq\varepsilon\|\varphi\|^2\,.
\end{multline*}

It follows that ${\|\Pi^{B_{\varepsilon}}\varphi\|}<\varepsilon {\|\varphi\|}$
and, since $g_n \to 1$ uniformly on $B$, also
%добавлено в англ.версию: поменяли $B'$ на $B$
that ${\|P^{(n)}\varphi\|}<\varepsilon{\|\varphi\|}$ if $n$
is sufficiently large. Since $\varepsilon$
is arbitrary, $$\|P^{(n)}\varphi\|\to 0\text{ as }n\to\infty,$$ and the proposition
is proved completely.

\section{Weak compactness of families of determinantal measures}
\subsection{Configurations and finite measures}

In a similar way as the Bessel point process of Tracy and Widom is the weak limit of its finite-dimensional approximations, the infinite determinantal measure ${\tilde {\mathbb B}}^{(s)}$, the sigma-finite analogue of the
Bessel point process for the values of $s$ smaller than $-1$, will be seen to be the scaling limit of its finite dimensional approximations, the infinite analogues of the Jacobi polynomial ensembles.
In this section, we  develop the formalism necessary for obtaining scaling limits of infinite determinantal measures. To do so, we will multiply our measures by finite densities, normalize and establish convergence of the resulting  determinantal
probability  measures. Proposition~\ref{optomeas} tells us that for finite determinantal measures induced by projection operators, local trace class convergence of the operators implies weak convergence of the determinantal measures (considered as measures on the space of Radon measures on the phase space). In order to prove the vanishing of the ``Gaussian parameter'' and to establish convergence of finite-dimensional approximations on the Pickrell set, we will however need a finer notion of convergence of probability measures on spaces of configurations: namely, under some additional assumptions we will code configurations by finite measures and determinantal measures by measures on the space of probability measures on the phase space. We proceed to precise definitions.

Let $f$ be a nonnegative measurable function on $E$, set
$$
\Conf_f(E)=\{X: \sum\limits_{x\in X}f(x)<\infty\},
$$
and introduce a map
 $\sigma_f:\Conf_f(E)\to \Mfin(E)$
by the formula $$
\sigma_f(X)=\sum\limits_{x\in X}f(x)\delta_x.
$$
(where $\delta_x$ stands, of course, for the
delta-measure at $x$).

Recall that the {\it intensity} $\xi\Prob$ of a probability measure $\Prob$ on
$\Conf(E)$ is a sigma-finite measure on $E$ defined, for a bounded Borel set $B\subset E$, by the formula
$$
\xi\Prob(B)=\int\limits_{\Conf(E)}\#_{B}(X)d\Prob(X).
$$
In particular, for a determinantal measure $\Prob_K$ corresponding to an
operator $K$ on $L_2(E, \mu)$ admitting a continuous kernel $K(x,y)$, the intensity
is, by definition,  given by the formula
$$
\xi\Prob_K=K(x,x)\mu.
$$
By definition, we have the following
\begin{proposition}\label{f-int}
Let $f$ be a nonnegative continuous function on $E$, and let $\Prob$ be a probability measure  on
$\Conf(E)$. If  $f\in L_1(E, \xi\Prob)$, then $\Prob(\Conf_f(E))=1$.
\end{proposition}
Under the assumptions of Proposition \ref{f-int}, the map $\sigma_f$ is $\Prob$-almost surely well-defined on $\Conf(E)$, and the
measure  $(\sigma_f)_*\Prob$ is a Borel probability measure on the space
$\Mfin\left(E\right)$, that is, an element of the space
$\Mfin\left(\Mfin(E)\right)$.

\subsection{Weak compactness and weak convergence in the space of configurations and in the space of finite measures}

We start by formulating a tightness criterion for such families of measures.
\begin{proposition}\label{tight-mfin}
Let $f$ be a nonnegative continuous function on $E$.
Let $\{\Prob_{\alpha}\}$ be a family of Borel probability measures on $\Conf(E)$ such that
\begin{enumerate}
\item $f\in L_1(E, \xi\Prob_{\alpha})$ for all $\alpha$ and
$$
\sup\limits_{\alpha} \int\limits_E fd\xi\Prob_{\alpha}<+\infty;
$$
\item  for any
$\varepsilon>0$ there exists a compact set $B_{\varepsilon}\subset
E$ such that $$\sup\limits_{\alpha}  \int\limits_{E\setminus B_{\varepsilon}} fd\xi\Prob_{\alpha}<\varepsilon.$$
\end{enumerate}

Then the family $(\sigma_f)_*\Prob_{\alpha}$ is tight in $\Mfin\left(\Mfin(E)\right)$.
\end{proposition}
%Our next aim is to give a sufficient condition for precompactness of families of such
{\bf Remark.} The assumptions of Proposition \ref{tight-mfin} can be equivalently reformulated as follows:
 the measures $(\sigma_f)_*\Prob_{\alpha}$ are all well-defined and the  family $f\xi\Prob_{\alpha}$ is tight in $\Mfin(E)$.

Proof of Proposition \ref{tight-mfin}. Given $\varepsilon>0$, our aim is to find a compact set
 $C\subset \Mfin(E)$ such that $(\sigma_f)_*\Prob_{\alpha}(C)>1-\varepsilon$ for all $\alpha$.

Let $\varphi:E\to\mathbb{R}$ be a bounded  function.
Define a measurable function
$\Int_{\varphi}:\Mfin(E)\to\mathbb{R}$ by the formula
$$\Int_{\varphi}(\eta)=\int_E\varphi d\eta.$$
Given a Borel subset $A\subset E$, for brevity we write
$\Int_A(\eta) = \Int_{\chi_A}(\eta)$.
%добавлено в англ.версию: (\eta)

The following proposition is immediate from local compactness of the space $E$ and the weak compactness
of the space of Borel probability measures on a compact metric space.
\begin{proposition}Let $L>0,\ \varepsilon_n>0,\
\lim\limits_{n\to\infty}\varepsilon_n=0.$ Let $K_n\subset E$ be
compact sets such that $\bigcup\limits_{n=1}^{\infty}K_n=E$. The
set $$\{\eta\in\Mfin(E):\Int_E(\eta)\leq
L, \Int_{E\backslash K_n}(\eta)\leq \varepsilon_n\ for\ all\
n\in\mathbb{N}\}$$ is compact in the weak topology on
$\Mfin(E).$
\end{proposition}
The Prohorov Theorem together with the Chebyshev Inequality now immediately implies
\begin{corollary}\label{cpct-set-mfin} Let $L>0,\ \varepsilon_n>0,
\ \lim\limits_{n\to\infty}\varepsilon_n=0.$ Let $K_n\subset E$ be
compact sets such that $\bigcup\limits_{n=1}^{\infty}K_n=E.$ Then
the set
\begin{multline*}
\Biggl\{\nu\in\Mfin(\Mfin(E)):\
\int\limits_{\Mfin(E)}\Int_E(\eta)d\nu(\eta)\leq L,\\
\int\limits_{\Mfin(E)}\Int_{E\backslash K_n}(\eta)d\nu(\eta)
\leq \varepsilon_n\text{ for all }n\in\mathbb{N}\Biggr\}
\end{multline*}
is compact in the weak topology on $\Mfin(\Mfin(E))$.
%добавлено в англ.версию: в пред.строке поменяли \Mfin(E) на \Mfin(\Mfin(E))
\end{corollary}

Corollary \ref{cpct-set-mfin} implies Proposition \ref{tight-mfin}.
First, the total mass of the measures $f\xi\Prob_{\alpha}$ is uniformly bounded, which, by the Chebyshev inequality,
 implies, for any $\varepsilon>0$,  the existence of the constant $L$ such that for all $\alpha$ we have
$$
(\sigma_f)_*\Prob_{\alpha}(\{\eta\in \Mfin(E): \eta(E)\leq L\})>1-\varepsilon.
$$

 Second, tightness of the family
$f\xi\Prob_{\alpha}$ precisely gives, for any $\varepsilon>0$, a compact set $K_{\varepsilon}\subset E$
satisfying, for all $\alpha$,  the inequality
$$
\int\limits_{\Mfin(E)}\Int_{E\backslash K_{\varepsilon}}(\eta)d(\sigma_f)_*\Prob_{\alpha}(\eta)
\leq \varepsilon.
$$
Finally, choosing a sequence $\varepsilon_n$ decaying fast enough and using Corollary \ref{cpct-set-mfin}, we conclude
the proof of Proposition \ref{tight-mfin}.

We now give sufficient conditions ensuring that convergence in the space of measures on the space of configurations
implies convergence of corresponding measures on the space of finite measures.
\begin{proposition}\label{wconv-mfin}
Let $f$ be a nonnegative continuous function on $E$.
Let $\Prob_n, n\in{\mathbb N}$, $\Prob$ be Borel probability measures on $\Conf(E)$ such that
\begin{enumerate}
\item $\Prob_{n}\to\Prob$ with respect to the weak topology on $\Mfin(\Conf(E))$ as $n\to\infty$;
\item $f\in L_1(E, \xi\Prob_n)$ for all $n\in {\mathbb N}$;
\item the family $f\xi\Prob_n$ is a tight family of finite Borel measures on $E$.
\end{enumerate}

Then $\Prob(\Conf_f(E))=1$ and
the measures $(\sigma_f)_*\Prob_n$ converge to $(\sigma_f)_*\Prob$ weakly
in $\Mfin\left(\Mfin(E)\right)$
as $n\to\infty$.
\end{proposition}

Proposition \ref{wconv-mfin} easily follows from Proposition  \ref{tight-mfin}.
First, we restrict ourselves to the open subset $\{x\in E: f(x)>0\}$ which itself is a complete separable
metric space with respect to the induced topology.
Next observe that the total mass of the measures $f\xi\Prob_{n}$ is uniformly bounded, which, by the
%добавлено в англ.версию: поменяли в пред. строке на $f\xi\Prob_{n}$
Chebyshev inequality, implies, for any $\varepsilon>0$,  the existence of the constant $L$ such that for all $n$ we have
$$
\Prob_{n}\left(\{X\in\Conf(E): \sum\limits_{x\in X}f(X)\leq L\}\right)>1-\varepsilon.
$$
Since the measures $\Prob_n$ converge to $\Prob$ weakly in $\Mfin(\Conf(E))$ and the set
$\{X\in\Conf(E): \sum\limits_{x\in X}f(X)\leq L\}$ is closed in $\Conf(E)$, it follows that
$$
\Prob\left(\{X\in\Conf(E): \sum\limits_{x\in X}f(X)\leq L\}\right)>1-\varepsilon,
$$
and, consequently, that $\Prob(\Conf_f(E))=1$, and the measure $(\sigma_f)_*\Prob$ is well-defined.

The family $(\sigma_f)_*\Prob_n$
is tight and must have a weak accumulation point $\Prob'$.

Using the weak convergence $\Prob_n\to \Prob$ in $\Mfin(\Conf(E))$, we now show that
the finite-dimensional distributions of $\Prob'$  coincide with those of $(\sigma_f)_*\Prob$.
Here we use the assumption that our function $f$ is positive and, consequently,
bounded away from zero on every bounded subset of our locally compact space $E$.

Indeed, let $\varphi_1, \dots, \varphi_l:E \to {\mathbb R}$ be  continuous  functions  with disjoint compact supports.

By definition, the joint distribution of the random variables
$\Int_{\varphi_1}, \dots, \Int_{\varphi_l}$ with respect to $(\sigma_f)_*\Prob_n$ coincides
with the joint distribution of the random variables $\#_{\varphi_1/f}, \dots, \#_{\varphi_l/f}$
with respect to $\Prob_n$. As $n\to\infty$, this joint distribution converges to
the joint distribution of  $\#_{\varphi_1/f}, \dots, \#_{\varphi_l/f}$
with respect to $\Prob$ which on the one hand, coincides with the the joint distribution of the random variables
$\Int_{\varphi_1}, \dots, \Int_{\varphi_l}$ with respect to $(\sigma_f)_*\Prob$ and, on the other hand, also coincides with the joint distribution of the random variables
$\Int_{\varphi_1}, \dots, \Int_{\varphi_l}$ with respect to $\Prob'$.

By Proposition \ref{fin-distr-unique},  the finite-dimensional distributions determine a measure uniquely.  Therefore,
$$\Prob'=(\sigma_f)_*\Prob, $$
and the proof is complete.

\subsection{Applications to determinantal point processes}

Let $f$ be a nonnegative  continuous function on $E$.
If an operator $K\in \scrI_{1,\loc}(E,\mu)$ induces
a determinantal measure $\Prob_K$ and satisfies
$fK \in \scrI_{1}(E,\mu)$, then
\begin{equation}\label{sup-f}
\Prob_K(\Conf_f(E))=1.
 \end{equation}
If, additionally, $K$ is assumed to be self-adjoint, then the
weaker requirement $\sqrt{f}K\sqrt{f} \in \scrI_{1}(E,\mu)$
also implies  (\ref{sup-f}).

In this special case, a sufficient condition for tightness takes the following form.

\begin{proposition}\label{det-prec} Let $f$ be a bounded nonnegative  continuous function on $E$.
Let $K_{\alpha}\in\scrI_{1,\loc}(E,\mu)$ be a family of self-adjoint positive contractions  such that
$$
\sup\limits_{\alpha} \tr \sqrt{f}K_{\alpha}\sqrt{f}<+\infty
$$
and such that for any
$\varepsilon>0$ there exists a bounded set $B_{\varepsilon}\subset
E$ such that $$\sup\limits_{\alpha}\tr\chi_{E\backslash
B_{\varepsilon}}\sqrt{f}K_{\alpha}\sqrt{f} \chi_{E\backslash
B_{\varepsilon}} <\varepsilon.$$
Then the  family of measures $\left\{\left(\sigma_{f}\right)_*\Prob_{K_{\alpha}}\right\}$ is weakly precompact in $\Mfin\left(\Mfin(E)\right)$.
\end{proposition}

\subsection{Induced processes corresponding to functions assuming values in $[0,1]$}
Let $g:E\to [0,1]$
be a nonnegative Borel function, and, as before, let $\Pi\in \scrI_{1, \loc}(E, \mu)$  be an orthogonal projection operator with range $L$ inducing a determinantal measure $\Prob_{\Pi}$ on $\Conf(E)$. Since the values of $g$ do not exceed $1$, the multiplicative functional $\Psi_g$ is automatically integrable. In this particular case Proposition \ref{pr1-bis} can be
reformulated as follows:
\begin{proposition} \label{pr1-bis-one} If
$\sqrt{1-g}\Pi\sqrt{1-g}\in \scrI_1(E, \mu)$   and
$||{(1-g)}\Pi||<1$, then
\begin{enumerate}
\item $\Psi_g$ is positive on a set of positive measure;
\item the subspace $\sqrt{g}L$ is closed, and the operator $\Pi^g$ of orthogonal projection onto the subspace $\sqrt{g}L$ is locally of trace class;
\item  we have
\begin{equation*}
\Prob_{\Pi^g}=\frac{\displaystyle \Psi_g\Prob_{\Pi}}
{\displaystyle \int\limits_{\Conf(E)}\Psi_g\,d\Prob_{\Pi}}.
\end{equation*}
\end{enumerate}
\end{proposition}
{\bf Remark.} Since the operator $\sqrt{1-g}\Pi$ is, by assumption,  Hilbert-Schmidt, and the the values of $g$ do not exceed $1$, the condition $||{(1-g)}\Pi||<1$ is equivalent to the condition  $||\sqrt{1-g}\Pi||<1$ and both are equivalent to the nonexistence
of a function $\Phi\in L$ supported on the set $\{x\in E: g(x)=1\}$. In particular, if the function $g$ is strictly positive and strictly less than 1, the condition is automatically verified.
Proposition \ref{det-prec} now implies
\begin{corollary}\label{cor1-bis-one}
Let $f$ be a bounded nonnegative continuous function on $E$. Under the assumptions of Proposition \ref{pr1-bis-one}, if
$$
\tr \sqrt{f} \Pi\sqrt{f} <+\infty,
$$
then also
$$
\tr \sqrt{f} \Pi^g\sqrt{f} <+\infty.
$$
\end{corollary}
Proof: Equivalently, we must prove that if the operator $\sqrt{f}\Pi$ is Hilbert-Schmidt, then the operator $\sqrt{f}\Pi^g$ is also Hilbert-Schmidt.
Since $\Pi^g=\sqrt{g}\Pi(1+(g-1)\Pi)^{-1}\sqrt{g}$, the statement is immediate from the fact that Hilbert-Schmidt operators form an ideal.

\subsection{Tightness for families of induced processes}
We now give a sufficient condition for the tightness of families of measures of the form $\Pi^g$ for fixed $g$.
This condition will subsequently be used for establishing convergence of determinantal measures obtained as products of infinite determinantal measures and multiplicative functionals.

Let $\Pi_{\alpha}\in \scrI_{1, \loc}(E, \mu)$ be a family of orthogonal projection operators
in $L_2(E, \mu)$.  Let $L_{\alpha}$ be the range of $\Pi_{\alpha}$. Let $g:E\to [0,1]$ be a Borel function such that for each $\alpha$
the assumptions of Proposition \ref{pr1-bis-one} are satisfied and thus the operators $\Pi_{\alpha}^g$ and the corresponding determinantal
measures   $\Prob_{\Pi_{\alpha}^g}$ are well-defined for all $\alpha$.
Furthermore,  let $f$ be a nonnegative function on $E$ such that
 such that for all $\alpha$ we have
\begin{equation}\label{tr-lim-alpha}
\sup\limits_{\alpha} \tr \sqrt{f}\Pi_{\alpha}\sqrt{f}<+\infty
\end{equation}
and such that for any
$\varepsilon>0$ there exists a bounded Borel set $B_{\varepsilon}\subset
E$ such that
\begin{equation}\label{tr-eps-alpha}
\sup\limits_{\alpha}\tr\chi_{E\backslash
B_{\varepsilon}}\sqrt{f}\Pi_{\alpha}\sqrt{f} \chi_{E\backslash
B_{\varepsilon}} <\varepsilon.
\end{equation}
(in other words, $f$ is such that all the assumptions of
Proposition \ref{det-prec} are satisfied for all $\alpha$). It follows from Corollary \ref{cor1-bis-one} that the measures $(\sigma_f)_*\Prob_{\Pi_{\alpha}^g}$ are also well-defined for all $\alpha$.

Sufficient conditions for tightness of this family of operators are given in the following

\begin{proposition} \label{g-det-prec}
In addition to the requirements, for all $\alpha$, of Proposition \ref{det-prec} and Proposition \ref{pr1-bis-one}, make the assumption
\begin{equation} \label{unif-ind-2}
\inf\limits_{\alpha} 1-|| (1-g)\Pi_{\alpha} ||>0.
\end{equation}
Then the family of measures $\left\{\left(\sigma_f\right)_*\Prob_{\Pi_{\alpha}^g}\right\}$ is weakly precompact in $\Mfin\left(\Mfin(E)\right)$.
\end{proposition}
Proof. The requirement (\ref{unif-ind-2}) implies that the norms of the operators $$(1+(g-1)\Pi_{\alpha})^{-1}$$ are uniformly bounded in $\alpha$.
Recalling that $\Pi_{\alpha}^g=\sqrt{g}\Pi_{\alpha}(1+(g-1)\Pi_{\alpha})^{-1}\sqrt{g}$, we obtain that (\ref{tr-lim-alpha}) implies
\begin{equation*}
\sup\limits_{\alpha} \tr \sqrt{f}\Pi_{\alpha}^g\sqrt{f}<+\infty,
\end{equation*}
while (\ref{tr-eps-alpha}) implies
\begin{equation*}
\sup\limits_{\alpha}\tr\chi_{E\backslash
B_{\varepsilon}}\sqrt{f}\Pi_{\alpha}^g\sqrt{f} \chi_{E\backslash
B_{\varepsilon}} <\varepsilon.
\end{equation*}

Proposition \ref{g-det-prec}  is  now immediate from Proposition \ref{det-prec}.

\subsection{Tightness of families of finite-rank deformations}

We next remark that, under certain additional assumptions,  tightness is preserved by
taking finite-dimensional deformations of determinantal processes.

As before, we let $\Pi_{\alpha}\in \scrI_{1, \loc}(E, \mu)$ be a family of orthogonal projection operators
in $L_2(E, \mu)$.  Let $L_{\alpha}$ be the range of $\Pi_{\alpha}$.
Let $v_{(\alpha)}\in L_2(E, \mu)$ be orthogonal to $L_{\alpha}$,  let $L_{\alpha}^{v}=L_{\alpha}\oplus {\mathbb C}v_{(\alpha)}$, and let
$\Pi_{\alpha}^{v}$ be the corresponding orthogonal projection operator. By the Macch{\`\i}-Soshnikov theorem, the operator  $\Pi_{\alpha}^{v}$ induces a determinantal measure $\Prob_{\Pi_{\alpha}^{v}}$ on $\Conf(E)$.
As above, we require that all the assumptions of
Proposition \ref{det-prec} be satisfied for the family $\Pi_{\alpha}$. The following Corollary is immediate from Proposition \ref{det-prec}.
\begin{proposition}\label{rank-one-orth}
Assume additionally
that the family of measures $f|v^{(\alpha)}|^2\mu$ is precompact in $\Mfin(E)$. Then the family of measures $\left\{\left(\sigma_f\right)_*\Prob_{\Pi_{\alpha}^v},\right\}$ is weakly precompact in $\Mfin\left(\Mfin(E)\right)$.
\end{proposition}
This proposition can be extended to perturbations of higher rank.  The assumption of orthogonality of $v^{\alpha}$ to $L_{\alpha}$ is too restrictive and can be weakened to an assumption that the angle between the vector and the subspace is bounded below: indeed, in that case we can orthogonalize and apply Proposition \ref{rank-one-orth}.

We thus take $m\in {\mathbb N}$ and assume that,
in addition to the family of $\Pi_{\alpha}$ of locally trace-class projection operators considered above, for every $\alpha$  we are given vectors $v^{(1)}_{\alpha},
\dots, v^{(m)}_{\alpha}$ of unit length, linearly independent and independent from $L_{\alpha}$. Set $$L_{\alpha}^{v, m}=L_{\alpha}\oplus {\mathbb C}v^{(1)}_{\alpha}\oplus\dots\oplus{\mathbb C} v^{(m)}_{\alpha},$$
and let  $\Pi_{\alpha}^{v, m}$ be the corresponding projection operator.

By the Macch{\`\i}-Soshnikov theorem, the operator  $\Pi_{\alpha}^{v, m}$ induces a determinantal measure $\Prob_{\Pi_{\alpha}^{v, m}}$ on $\Conf(E)$.
As above, we require that all the assumptions of
Proposition \ref{det-prec} be satisfied for the family $\Pi_{\alpha}$. Here $\ang(v, L)$ stand for the angle between a nonzero vector $v$ and
a closed subspace $L$.
\begin{proposition} \label{fin-rank-general}
Assume additionally that
\begin{enumerate}
\item the family of measures $f|v_{\alpha}^{(k)}|^2\mu$, over all $\alpha$ and $k$, is precompact in $\Mfin(E)$;
\item there exists $\delta>0$ such that for any $k=1, \dots, m$  and all $\alpha$ we have
$$
\ang(v_{\alpha}^{(k)}, L_{\alpha}\oplus {\mathbb C}v^{(1)}_{\alpha}\oplus\dots\oplus{\mathbb C} v^{(k-1)}_{\alpha})\geq \delta.
$$
\end{enumerate}
Then the family of measures $\left\{(\sigma_f)_*\Prob_{\Pi_{\alpha}^{v, m}},\right\}$ is weakly precompact in $\Mfin\left(\Mfin(E)\right)$.
\end{proposition}

The proof proceeds by induction on $m$. For $m=1$, it suffices to apply Proposition \ref{rank-one-orth}
to the vector obtained by taking the orthogonal projection of $v_{\alpha}^{(1)}$ onto the orthogonal complement of $L$. For the induction step, similarly, we apply Proposition \ref{rank-one-orth} to the vector obtained by taking the orthogonal projection of $v_{\alpha}^{(m)}$ onto the orthogonal complement of
$L_{\alpha}\oplus {\mathbb C}v^{(1)}_{\alpha}\oplus\dots\oplus{\mathbb C} v^{(m-1)}_{\alpha}$. The proposition is proved completely.

\subsection{Convergence of finite-rank perturbations}

A sufficient condition for weak convergence of determinatal measures considered as elements of the space $\Mfin(\Mfin(E))$ can be formulated as follows.

\begin{proposition}\label{weak-conv-f} Let $f$ be a nonnegative continuous function on $E$.
Let $K_n, K\in \scrI_{1,\loc}$ be self-adjoint positive contractions such that
$K_n\to K$ in
$\scrI_{1,\loc}(E, \mu)$ as $n\to\infty$.
Assume additionally that
\begin{equation}\label{f-trace}
\sqrt{f}K_n\sqrt{f}\to \sqrt{f}K\sqrt{f} \ \text{in} \ \scrI_{1}(E, \mu)
\end{equation}
as $n\to\infty$.
Then
$$(\sigma_f)_*\Prob_{K_n}\to(\sigma_f)_*\Prob_{K}$$
weakly in $\Mfin(\Mfin(E))$ as $n\to\infty$.
\end{proposition}
%добавлено в англ.версию: меняем в пред. строке $\mathfrak{M}(\Mfin(E))$ на $\Mfin(\Mfin(E))$

Combining Proposition \ref{weak-conv-f} with, on the one hand, Propositions \ref{g-det-prec}, \ref{fin-rank-general} and, on the other hand,  Propositions \ref{conv-induced}, \ref{fin-pert} and \ref{conv-inf-det}, we arrive at the following
\begin{proposition}\label{conv-unif-finrank}
\begin{enumerate}

\item In the notation and under the assumptions of Proposition \ref{conv-induced},  additionally  require (\ref{f-trace}) to hold.
Then we have $$\sqrt{f}\tilde{\buf}(g, \Pi_n)\sqrt{f}\to \sqrt{f}\tilde{\buf}(g, \Pi)\sqrt{f}$$ in
$\scrI_{1}(E,\mu)$,
and, consequently,
$$
(\sigma_f)_*\Prob_{\tilde{\buf}(g, \Pi_n)}\to(\sigma_f)_*\Prob_{\tilde{\buf}(g, \Pi)}
$$
with respect to the weak topology on $\Mfin(\Mfin(E))$ as $n\to\infty$.

\item In the notation and under the assumptions of Proposition \ref{fin-pert},  additionally  require (\ref{f-trace}) to hold.
Then we have
$$
\sqrt{f}\widetilde{\Pi}_n\sqrt{f}\to \sqrt{f}\widetilde{\Pi}\sqrt{f} \ \text{in} \ \scrI_{1}(E,\mu) \ \text{as}  \ n\to\infty,
$$
and, consequently, $(\sigma_f)_*\Prob_{\widetilde{\Pi}_n} \to (\sigma_f)_*\Prob_{\widetilde{\Pi}}$ with respect to the weak topology on  $\Mfin(\Mfin(E))$ as $n\to\infty$;
%добавлено в англ.версию: поменяно выше $\Prob_{\widetilde{\Pi}_n}\to\Prob_{\widetilde{\Pi}}$ на
% $(\sigma_f)_*\Prob_{\widetilde{\Pi}_n} \to (\sigma_f)_*\Prob_{\widetilde{\Pi}}$

\item In the notation and under the assumptions of Proposition \ref{conv-inf-det},  additionally  require (\ref{f-trace}) to hold.
Then we have
$$
\sqrt{f}\Pi^{(g,n)}\sqrt{f}\to\sqrt{f}\Pi^g \sqrt{f}\;\text{ in }\; \scrI_{1}(E,\mu)
\;\text{ as }\; n\to\infty\;.
$$
and, consequently,
$$
(\sigma_f)_*\frac{\Psi_g\mathbb{B}^{(n)}}{\displaystyle \int\Psi_g\,d\mathbb{B}^{(n)}} \to (\sigma_f)_*\frac{\Psi_g\mathbb{B}}{\displaystyle \int\Psi_g\,d\mathbb{B}}
$$ with respect to the weak topology on  $\Mfin(\Mfin(E))$ as $n\to\infty$.
\end{enumerate}
\end{proposition}

\appendix

{\bf Acknowledgements.} This project has received funding from the European Research Council (ERC) under the European Union's Horizon 2020 research and innovation programme (grant agreement No 647133 (ICHAOS))
and has also been funded by the Russian Academic Excellence Project `5-100'.

\end{document}